\documentclass[10pt,a4paper]{article}
\usepackage{color}
\usepackage{ifpdf}
\ifpdf
 \usepackage[pdftex]{graphicx}
\else
 \usepackage[dvips]{graphicx}
\fi

\usepackage{hyperref}

\usepackage[T1]{fontenc}
\usepackage{a4,amssymb,amsthm,amsmath}
\usepackage{lmodern}
\usepackage[all]{xy}
\def\à{\`a}
\def\è{\`e}
\def\ä{\"a}
\newcommand{\R}{\mathbb{R}}
\newcommand{\C}{\mathbb{C}}
\newcommand{\T}{\mathbb{T}}
\newcommand{\PP}{\mathbb{P}}

\newcommand{\Sp}{\mathbb{S}}
\newcommand{\Z}{\mathbb{Z}}
\newcommand{\N}{\mathbb{N}}
\newcommand{\lra}{\longrightarrow}
\newcommand{\lms}{\longmapsto}
\newcommand{\bw}{\bigwedge}
\newcommand{\w}{\wedge}
\theoremstyle{definition}

\DeclareMathOperator{\Vect}{Vect}
\DeclareMathOperator{\Kod}{Kod}
\DeclareMathOperator{\End}{End}
\DeclareMathOperator{\Aut}{Aut}

\author{Guillaume~Deschamps} \title{Compatible Complex
  Structures on Twistor Spaces} \date{\today}

\begin{document}

\maketitle
\begin{center}
guillaume.deschamps@univ-brest.fr
\end{center}
\quad\\
{\it Let $(M,g)$ be a Riemannian
  $4$-manifold.
  The twistor space $Z\to M$ is a $\C P^1$-bundle
  whose total space $Z$ admits a natural metric $\tilde
  g$. The aim of this article is to study properties of
  complex structures on $(Z,\tilde g)$ which are compatible
  with the $\C P^1$-fibration and the metric $\tilde g$. The
  results obtained enable us to translate some metric
  properties on $M$ in terms of complex properties on its
  twistor space $Z$.}

\section*{Introduction}
Let $(M,g)$ be an oriented 4-dimensional Riemannian manifold (not
necessarily compact). Due to the Hodge-star operator $\star$, we
have a decomposition of the bivector bundle
$\bigwedge^2TM=\bigwedge^+\oplus\bigwedge^-$. Here $\bigwedge^\pm$
is the eigen-subbundle for the eigenvalue $\pm 1$ of $\star$. The
metric $g$ on $M$ induces a metric, denoted by
$<\phantom{.},\phantom{.}>$, on the bundle $\bigwedge^2 TM$. Let
$\pi:\,Z=\Sp\bigl(\bigwedge^+\bigr)\lra M$ be the sphere bundle;
the fiber over a point $m\in M$ parameterizes the complex
structures on the tangent space $T_mM$ compatible with the
orientation and the metric $g$. It is the twistor space of the
manifold $(M,g)$. Since the structural group of $Z$ is $SO(3)\subset \Aut(\C P^1)$, we
can thus put the complex structure of $\C P^1$ on each fiber. On
the other hand, the Levi-Civita connection on $(M,g)$ induces a
splitting of the tangent bundle $TZ$ into the direct sum of the
horizontal and vertical distributions: $TZ=H\oplus V$.  Therefore,
the twistor space $Z$ admits a natural metric $\tilde g$ defined
by its restrictions to $H$ and $V$: we endow $V$ with the
Fubini-Study metric and $H\simeq \pi^\star TM$ with the pullback
of the metric $g$.

In this article we study some aspects of almost complex structures
on $(Z,\tilde g)$ which are Hermitian and extend the complex
structure of the fibers. These structures will be called {\it
compatible almost
  complex structures} on $(Z,\tilde g)$. In particular, the
integrability of two such structures means that the metric
$\tilde g$ is bihermitian \cite{Pon97}, \cite{AGG99}.

To each morphism respecting the twistor fibration
\[
\xymatrix @R=0,2pc @C=0,5cm{
  Z\ar[ddddr]_{\pi}&\ar[rr]^{f}&&& Z\ar[ddddl]_{\pi}\\
  \\
  \\
  \\
  && M&&\\
}
\]
we associate a compatible almost complex structure $\mathbb
J_f$ on $(Z,\tilde g)$ in the following way. Let $z\in Z$
with $\pi(z)=m\in M$, and write $T_zZ=H_z\oplus V_z$. Here,
$V_z$ is the tangent space to the fiber $\pi^{-1}(m)\simeq\C
P^1$ and is therefore equipped with a complex structure. On
the other hand, we endow $H_z\simeq T_m M$ with the complex
structure associated to the point $f(z)$. Conversely, any
compatible almost complex structure $\mathbb J$ on
$(Z,\tilde g)$ defines a unique morphism $f : Z\lra Z$
respecting the fibration such that $\mathbb J_f=\mathbb J$.

The almost complex structure $\mathbb J_{Id}$ associated to the
identity is the canonical twistor almost complex structure
\cite{AHS78}. If $\sigma$ is the morphism of $Z$ whose restriction
to each fiber of $\pi$ is the antipodal map of $\Sp^2$, we denote
by $\mathbb J_{\sigma}$ the almost complex structure associated to
$\sigma$. Now, an almost complex manifold $(M,g,J_M)$ such that
$J_M$ is compatible with the orientation and the metric $g$
defines a tautological section of $Z\lra M$. This section can be
taken as the infinity section and we can therefore consider the
constant morphism $f=\infty$. The associated almost complex
structure will be denoted by $\mathbb J_\infty$. Let
$\lambda\in\C^\star$ and consider the morphism $f=\lambda Id$
acting as $\lambda Id$ in each fiber minus infinity (i.e.~$\C
P^1-\{\infty\}\simeq\C$) and preserving infinity.  We  denote by
$\mathbb J_{\lambda Id}$ the corresponding almost complex
structure on $Z$.

The integrability of the structures $\mathbb J_{Id}$,
$\mathbb J_\sigma$, $\mathbb J_{\infty}$, $\mathbb
J_{\lambda Id}$ are related to the curvature of the
metric $g$ on $M$. Let $R : \bigwedge^2TM\lra\bigwedge^2TM$ be
the curvature operator. The decomposition
$\bigwedge^2TM = \bigwedge^+\oplus\bigwedge^-$ allows us to
write $R$ in block matrix form as follows
\[
R=\left(\begin{array}{cc} A&^tB\\ B&C\end{array}\right),
\]
where $A=W^++\frac{s}{12}Id$, $C=W^-+\frac{s}{12}Id$, $W^+$
({\it resp.}~$W^-$) is the selfdual ({\it
  resp.}~anti-selfdual) Weyl tensor, $s$ is the scalar
curvature and $B$ the trace-free Ricci curvature
\cite{Bes87}.

\smallskip

The main result of this article is the following:

\medskip

\noindent {\bf Theorem 1.}  \textit{Let $(M,g)$ be an
  oriented Riemannian
  4-manifold.
\begin{enumerate}
\item[{\rm A)}] The complex structure $\mathbb J_\sigma$ is never
  integrable.
\item[{\rm B)}] The complex structure $\mathbb J_{Id}$ is
  integrable if, and only if, $g$ anti-selfdual (i.e. $A$ is a homothety)
  {\em \cite{AHS78}}.
\item[{\rm C)}] Let $J_M$ be an almost complex structure on $M$
  compatible with the metric $g$ and the orientation. The
  complex structure $\mathbb J_{\infty}$ is integrable  if,
  and only if:
\begin{enumerate}
\item[\rm i)] $J_M$ is integrable;
\item[\rm ii)] the kernel of $A$ contains the plane
  $J_M^\perp\subset \bigwedge^+$ orthogonal to the line
  generated by $J_M$.
\end{enumerate}
\item[{\rm D)}] Let $(M,g,J_M)$ be a K\"ahlerian manifold.
If $\lambda \notin\{0,1\}$, the complex
  structure $\mathbb J_{\lambda Id}$ is integrable if, and only if,
  $(M,g,J_M)$ is scalar-flat Kähler (i.e.~$A$=0).
\item[{\rm E)}] Let $(M,g)$ be an anti-selfdual Riemannian
  manifold. Its scalar curvature is zero if, and only if,
  any $m\in M$ has an open neighborhood $\mathcal U$ such
  that, over $\mathcal{U}$, $(Z,\tilde g)$ admits a
  compatible complex structure
different from $\mathbb J_{Id}$.
\end{enumerate}
}

\medskip

The conditions i)~\&~ii) of part~C in the previous theorem are
satisfied as soon as $(M,g,J_M)$ is Kähler. We show in section~C
that this K\"ahlerian property is equivalent
 to the integrability of
$\mathbb J_{\infty}$ in the compact case. For a scalar-flat
Kähler manifold $(M,g,J_M)$, the complex structures $\mathbb
J_{Id}$ \cite{Gau81}, $\mathbb J_\infty$ and $\mathbb
J_{\lambda Id}$ are integrable and compatible with the
metric $\tilde g$ on $Z$. This gives us a huge family of
real $6$-dimensional manifolds admitting a bihermitian
metric.

Recall that the Penrose correspondence gives a dictionary between
holymorphic properties of the twistor space $Z$ and properties of
the Riemannian manifold $(M,g)$. The above result can be viewed as
a new paragraph of that dictionary. In particular, we deduce from
it some new characterizations of K\"ahlerian metrics,
anti-selfdual scalar-flat metrics and scalar-flat Kähler metrics,
in terms of twistor spaces.

The proof of Theorem~1 is split into five theorems,
Theorem~A,\dots, E, the proof of each being given in the
corresponding labelled section.

In section~D we explain how Theorem~1 can be used in order
to build a $1$-dimensional family of non conformal
anti-selfdual metrics on a scalar-flat Kähler manifold, or a
$1$-dimensional family of biholomorphisms on its
corresponding complex twistor space $(Z,\mathbb J_{Id})$.

In section~F we study more precisely the set of all compatible
complex structures on the twistor space  of a locally conformally
K\"ahler manifolds. Whereas on section~G we will study the case of
bielliptic surfaces.

We conclude the paper by giving a generalisation of this theorem
to quaternionic K\"ahler manifolds of dimension $4n$ for $n>1$.

\subsection*{Notation}
We will use Einstein summation convention over repeated indices.
The fiber of $\pi : Z\lra M$ over $m\in M$ will be freely
identified with $\Sp^2$, $\C P^1$ or $SO(4)/U(2)$, the set of all
complex structure on $T_m M$. The bundle of bivectors $\bigwedge^2
TM$ will be identified with the bundle of skew-symmetric
endomorphisms of $TM$, or to the bundle of $2$-forms.

Let $(\theta_1,\theta_2,\theta_3,\theta_4)$ be an oriented
$g$-orthonormal frame defined over an open set $\mathcal U$
of $(M,g)$. Define three linear operators $I, J, K\in
\End(TM)$, over $\mathcal U$, by their matrix in the basis
$(\theta_1,\dots,\theta_4)$:
$$
I=\left[\begin{array}{cccc} 0&-1&0&0\\
    1&0&0&0\\
    0&0&0&-1\\
    0&0&1&0
  \end{array}
\right]\; J=\left[\begin{array}{cccc}
    0&0&-1&0\\
    0&0&0&1\\
    1&0&0&0\\
    0&-1&0&0
  \end{array}
\right]\; K=\left[\begin{array}{cccc}
    0&0&0&-1\\
    0&0&-1&0\\
    0&1&0&0\\
    1&0&0&0
  \end{array}
\right].
$$
Then, $(I,J,K)$ gives an oriented orthonormal basis over
$\mathcal U$ of $\bigwedge^+$ and therefore defines a
trivialization of the twistor space $\pi: Z\lra M$ over
$\mathcal U$:
$$
\pi^{-1}(\mathcal U)\simeq \mathcal U\times SO(4)/U(2).
$$
Let $(\theta^*_1,\dots,\theta^*_4)$ be the local coframe dual
to $(\theta_1,\dots,\theta_4)$.  Locally, the covariant
derivative $\nabla$ (on $M$) defined by the Levi-Civita
connection of the metric $g$ writes
$\nabla\theta_j=\Gamma_{ij}^k\theta^*_i\otimes\theta_k $.
The $\Gamma_{ij}^k$ are the Christoffel symbols of the
connection $\nabla$; they satisfy
$\Gamma_{ij}^k=-\Gamma_{ik}^j$.

Let $z\simeq(m,Q)\in\pi^{-1}(\mathcal U)$ be a point of $Z$
and write the tangent space as the direct sum of the
horizontal and vertical tangent spaces: $T_zZ=V_z\oplus
H_z$.  Denote by $\hat\theta\in H_z\simeq T_m M$ the
horizontal lift of $\theta\in T_m M$. We then have
\cite{dBN98}:
\[
\left\{
  \begin{array}{l}
    V_z=\left\{X\frac{\partial}{\partial
        Q} \, \mid \,  X\in \End(T_mM), \ ^tX=-X\textrm{ et }
      QX=-XQ\right\}\\
    H_z=\Vect\Big(\hat\theta_1(z),\dots,\hat\theta_4(z)\Big)
  \end{array}
\right.
\]
with $\left\{
  \begin{array}{l}
    \hat\theta_i(z)=\theta_i(m)-[\Gamma_{i\centerdot}^\centerdot(m),Q]
\frac{\partial}{\partial  Q}\\

    [\Gamma_{i\centerdot}^\centerdot(m),Q]\frac{\partial}{\partial Q}
    =\Big(\Gamma_{i\centerdot}^\centerdot(m)Q-Q\Gamma_{i\centerdot}^\centerdot(m)\Big)\frac{\partial}
{\partial Q}\in V_z.\end{array}\right.$

\quad\\
{\bf Remark}: The complex structure of rational curves on
the fiber $\pi^{-1}(m)\simeq\Sp^2$ at a point $z=(m,Q)$ is
given by the application \cite{dBN98}:
$$
\begin{array}{ccc}
  V_z\simeq T_Q\Sp^2&\lra&V_z\simeq T_Q\Sp^2\\
  X\frac{\partial}{\partial Q}&\lms&
  QX\frac{\partial}{\partial Q}.
\end{array}
$$

\quad\\
For all $A\in so(4)=\{A\in \End(TM) \, \mid \, ^tA=-A\}$ we can define
the vertical vector field $\tilde
A=[A,Q]\frac{\partial}{\partial Q}$. These vector fields
will be called {\it basic}.

\subsection*{A) General results}
In this section $(M,g)$ will be an oriented Riemannian
$4$-manifold. Results~-- and proofs --~given here in dimension
$4$, can be easily adapted to quaternionic K\"ahler $4n$-manifolds
and will be used in the last section of the paper.

To study the integrability of the almost complex structure
$\mathbb J_f$ we need to compute the Nijenhuis tensor $N$ of
$\mathbb J_f$ \cite{NN57}:
$$
N(X,Y)=[\mathbb J_f X,\mathbb J_fY]-\mathbb J_f[\mathbb
J_fX,Y]-\mathbb J_f[X,\mathbb J_fY]-[X,Y]\qquad\forall
(X,Y)\in T_z Z.
$$
The first necessary condition  for the integrability of $\mathbb
J_f$ appears in the next proposition.

\quad\\
{\bf Proposition 1.} {\it For any morphism $f$ we have:
\begin{enumerate}
  \item[{\rm i)}]  $N(X,Y)=0$ for all $X,Y\in V_z$;

 \item[{\rm ii)}] let $X,\theta\in V_z\times H_z$, then
\begin{itemize}
\item the vertical component of $N(X,\theta)$ is
        zero
\item the horizontal component of $N(X,\theta)$ is
        zero if and only if the restriction of
      $f$  to each fiber is holomorphic.
\end{itemize}
\end{enumerate}}
\quad\\
 As $\sigma$ is an anti-holomorphic involution on fibers we
easily get:

\quad\\
{\bf Theorem A.} {\it The almost complex structure
  $J_\sigma$ is never integrable.}

\quad\\
{\bf Proof of Proposition~1.} For any morphism $f$, each fiber of
$\pi : Z\lra M$ has  the structure of $\C P^1$. It follows
immediately from \cite{NN57} that $N(X,Y)=0$ for all $X,Y\in V_z$.

Let $\tilde X$ be a basic vertical vector field and
$\pi^{-1}(m)$ be a fixed fiber. The  restriction to that fiber
of the application $f$ is:
$$
\begin{array}{lccc}
  f\vert_{\pi^{-1}(m)}:&\Sp^2\simeq\pi^{-1}(m)&\lra&\Sp^2
  \simeq\pi^{-1}(m)\\  &Q&\lms& f(Q)
\end{array}
$$
Observe that $[\tilde X,\hat \theta_i]$ is vertical when
$\tilde X$ is. Since  the action of the complex structure
$\mathbb J_f$ on the fiber is equal to the rational curve
structure, it does not depend on the fiber. We then
have: $[\mathbb J_f\tilde X,\hat\theta_i]=[Q\tilde
X,\hat\theta_i]=Q[\tilde X,\hat\theta_i]=\mathbb J_f[\tilde
X,\hat\theta_i]$. This implies that, for
$i\in\{1,\dots,4\}$:
\[\begin{array}{lll}
  N(\tilde X,\hat\theta_i)&=&[Q\tilde X,f(Q)\hat\theta_i]-Q[Q\tilde
  X,\hat\theta_i]+\mathbb
  J_f[\tilde X,f(Q)\hat\theta_i]-[\tilde X,\hat\theta_i]\\
  &=&\Big((Q \tilde X).f(Q)\;\; -f(Q)(
  \tilde X.f(Q))\Big)\;\hat\theta_i\\
  &=&\Big(d_Qf(Q \tilde X)\;\;
  -f(Q) d_Qf(\tilde X)\Big)\;\hat\theta_i
\end{array}
\]
where $d_Qf$ is the differential of $f$ at $Q\in\Sp^2$.  The
horizontal component of $N(X,\theta)$ vanishes for all
$(X,\theta)\in V_z\times H_z$
if and only if the restrictions of $f$ to the fibers are
holomorphic .
\qed\\

In the trivialization of $Z\lra M$ over an open set
$\mathcal U$, the morphism $f$ can be written:
$$
\begin{array}{lccc}
  f\vert_{\pi^{-1}(\mathcal U)}:&\mathcal
  U\times\Sp^2&\lra&\mathcal U\times\Sp^2\\
  &(x,Q)&\lms& \Big(x,f(x,Q)\Big).
\end{array}
$$
In order to simplify the notation we set $P=f(x,Q)$ and
$[P_i^j]$ denotes the matrix, in the basis
$(\theta_1,\dots,\theta_4)$, of the operator $P$ viewed as
an endomorphism of $TM$.

\quad\\
{\bf Proposition 2.} {\it Let $f$ be any morphism and
  $(m,Q)\in Z$. Then, for all $i,j \in \{1,\dots,4\}$ one
  has:
\begin{enumerate}
 \item[{\rm i)}] the horizontal component of
  $N(\hat\theta_i,\hat\theta_j)$ can be written as
  $\widehat{E(\theta_i,\theta_j)}+F_{ij}$

  \item[{\rm ii)}] the vertical component of
  $N(\hat\theta_i,\hat\theta_j)$ can be written as $
  G(\theta_i,\theta_j)\frac{\partial}{\partial Q}$,
\end{enumerate}
  \quad\\
  where $\left\{\begin{array}{l} \;\;E(\theta_i,\theta_j)
      \textrm{
        is the  Nijenhuis tensor of the almost complex
        structure $P_0$}\\
      \quad\quad\textrm{on } TM \textrm{ defined  by }
      f(\centerdot,Q) \textrm{ over the
        open set }\mathcal U  \textrm{ (where  $Q$ is fixed);}\\
      \\
      \begin{array}{lll}
        F_{ij}&=&-P_i^r[\Gamma_{r\centerdot}^\centerdot,Q]\frac{\partial}
        {\partial
          Q}P_j^r\;\;\hat\theta_r
        +P_j^r[\Gamma_{r\centerdot}^\centerdot,Q]\frac{\partial}{\partial
          Q}P_i^l\;\;\hat\theta_l\\
        &&-P\Big([\Gamma_{j\centerdot}^\centerdot,Q]\frac{\partial}{\partial
          Q}P_i^l\;\;\hat\theta_l -[\Gamma_{i\centerdot}^\centerdot,Q]
\frac{\partial}{\partial Q}P_j^r\;\;\hat\theta_r\Big);
      \end{array}\\
      \\
      \;\;G(\theta_i,\theta_j)= \left[R\Big(\theta_i\wedge
        \theta_j-P\theta_i\wedge P\theta_j\Big)+
        QR\Big(P\theta_i\wedge\theta_j+\theta_i\wedge
        P\theta_j\Big),Q\right].
    \end{array}\right.
  $ }
\smallskip

\quad\\
{\bf Proof.}
 The curvature tensor is
$R(\theta_i,\theta_j)=\nabla_{\theta_i}\nabla_{\theta_j}-
\nabla_{\theta_j}\nabla_{\theta_i}-\nabla_{[\theta_i,\theta_j]}=
R_{kij}^l\theta_k^\star\otimes\theta_l$, with $
R_{kij}^l=g\Big(R(\theta_i,\theta_j)\theta_k,\theta_l\Big)$. Hence,
$$
R(\theta_i,\theta_j)\theta_k=\nabla_{\theta_i}(\Gamma_{jk}^m\theta_m)-
\nabla_{\theta_j}(\Gamma_{ik}^m\theta_m)-
\nabla_{(\Gamma_{ij}^m-\Gamma_{ji}^m)\theta_m}\theta_k
$$
yields
$$
R_{ijk}^l=\theta_i(\Gamma_{jk}^l)-\theta_j(\Gamma_{ik}^l)+[
\Gamma_{i\centerdot}^\centerdot,
\Gamma_{j\centerdot}^\centerdot]_k^l-(\Gamma_{ij}^\centerdot-\Gamma_{ji}^\centerdot)
\Gamma_{\centerdot k}^l.
$$
To finish the proof of the proposition we need the following
lemma.

\quad\\
{\bf Lemma 1.} {\it The Lie bracket of $\hat\theta_i$ with
  $\hat\theta_j$ satisfies:
$$[\hat\theta_i,\hat\theta_j]=\widehat{[\theta_i,\theta_j]}-
\left[R_{\centerdot
ij}^\centerdot,Q\right]\frac{\partial}{\partial Q}.$$

}
\quad\\
{\bf Proof of  Lemma~1.} From
$\hat\theta_i=\theta_i-[\Gamma_{i\centerdot}^\centerdot,Q]
\frac{\partial}{\partial
  Q}$ we can deduce that:
$$\begin{array}{lll}
  [\hat\theta_i,\hat\theta_j]&=&\left[\theta_i-[\Gamma_{i\centerdot}^\centerdot,Q]
\frac{\partial}{\partial
  Q}\;,\;\theta_j-[\Gamma_{j\centerdot}^\centerdot,Q]\frac{\partial}{\partial
  Q}\right]\\
  &=&[\theta_i,\theta_j]-[\theta_i(\Gamma_{j\centerdot}^\centerdot),Q]
\frac{\partial}{\partial
Q}+[\theta_j(\Gamma_{i\centerdot}^\centerdot),Q]
\frac{\partial}{\partial Q}
  -\Big[[\Gamma_{i\centerdot}^\centerdot,\Gamma_{j\centerdot}^\centerdot],
  Q\Big]\frac{\partial}{\partial
    Q}\\
  &=&\Big(\Gamma_{ij}^m-\Gamma_{ji}^m\Big)\theta_m-
\Big[[\theta_i(\Gamma_{j\centerdot}^\centerdot)-
  \theta_j(\Gamma_{i\centerdot}^\centerdot)+
  [\Gamma_{i\centerdot}^\centerdot,\Gamma_{j\centerdot}^\centerdot],Q\Big]
  \frac{\partial}{\partial
    Q}\\
  &=&\Big(\Gamma_{ij}^m-\Gamma_{ji}^m\Big)\theta_m-\Big(
  \Big[R_{\centerdot ij}^\centerdot,Q\Big]+(\Gamma_{ij}^m-\Gamma_{ji}^m)
  [\Gamma_{m\centerdot}^\centerdot,Q]
\Big)\frac{\partial}{\partial
    Q}\\
  &=&(\Gamma_{ij}^m-\Gamma_{ji}^m)\hat\theta_m-
  \left[R_{\centerdot ij}^\centerdot,Q\right]\frac{\partial}{\partial Q}\\
  &=&\widehat{[\theta_i,\theta_j]}-
  \left[R_{\centerdot ij}^\centerdot,Q\right]\frac{\partial}{\partial Q}. \quad
  \square
\end{array}
$$

\quad\\
We can now complete the proof of Proposition~1. The Nijenhuis
tensor is given by
$$\begin{array}{lll}
  N(\hat\theta_i,\hat\theta_j)&=&[\mathbb J_f \hat\theta_i,\mathbb
  J_f\hat\theta_j]-\mathbb J_f \Big([\mathbb
  J_f\hat\theta_i,\hat\theta_j]+[\hat\theta_i,\mathbb
  J_f\hat\theta_j]\Big)- [\hat\theta_i,\hat\theta_j],
\end{array}
$$
where:
$$\begin{array}{rll}
  [\mathbb J_f \hat\theta_i,\mathbb J_f\hat\theta_j]&=&
  [P_i^l\hat\theta_l,P_j^r\hat\theta_r]\\
  &=&
  \widehat{P\theta_i}.(P_j^r)\;\hat\theta_r-\widehat{P\theta_j}.(P_i^l)
\;\hat\theta_l
  +P_i^lP_j^r[\hat\theta_l,\hat\theta_r]\\
  \, [\mathbb J_f\hat\theta_i,\hat\theta_j]+[\hat\theta_i,\mathbb
  J_f\hat\theta_j]&=&[P_i^l\hat\theta_l,\hat\theta_j]+
[\hat\theta_i,P_j^r\hat\theta_r]\\
  &=&-\hat\theta_j.(
  P_i^l)\;\hat\theta_l+P_i^l[\hat\theta_l,\hat\theta_j]+
  \hat\theta_i.(
  P_j^r)\;\hat\theta_r+P_j^r[\hat\theta_i,\hat\theta_r].
\end{array}
$$
By Lemma~1 the horizontal component of the Nijenhuis tensor
is:
$$
\begin{array}{lll}
  \mathcal
  H\,N(\hat\theta_i,\hat\theta_j)&=&
\widehat{P\theta_i}.(P_j^r)\;\hat\theta_r
  -\widehat{P\theta_j}.(P_i^l)\;\hat\theta_l
  +P_i^lP_j^r\widehat{[\theta_l,\theta_r]}\\
  &&-P\Big(-\hat\theta_j.(
  P_i^l)\;\hat\theta_l+P_i^l\widehat{[\theta_l,\theta_j]}+

  \hat\theta_i.(
  P_j^r)\;\hat\theta_r+P_j^r\widehat{[\theta_r,\theta_i]}\Big)\;\;-\;\;
  \widehat{[\theta_i,\theta_j].}
\end{array}
$$
Fix $Q$ and denote by $P_0$ the almost complex structure on
$TM$, over $\mathcal{U}$, defined by $P_0(m)=f(m,Q)$. Then:
$$
\begin{array}{lll}
  \mathcal H\,N(\hat\theta_i,\hat\theta_j)
  &=&\widehat{[P_0\theta_i,P_0\theta_j]}-P_0
\Big(\widehat{[P_0\theta_i,\theta_j]}+
  \widehat{[\theta_i,P_0\theta_j]}\Big)-
\widehat{[\theta_i,\theta_j]}\\
  &&-P_i^r[\Gamma_{r\centerdot}^\centerdot,Q]\frac{\partial}{\partial
    Q}P_j^r\;\;\hat\theta_r
  +P_j^r[\Gamma_{r\centerdot}^\centerdot,Q]\frac{\partial}{\partial
    Q}P_i^l\;\;\hat\theta_l\\
  &&-P\Big([\Gamma_{j\centerdot}^\centerdot,Q]\frac{\partial}{\partial
    Q}P_i^l\;\;\hat\theta_l -[\Gamma_{i\centerdot}^\centerdot,Q]
\frac{\partial}{\partial Q}P_j^r\;\;\hat\theta_r\Big)\\
  &=&\widehat{E(\theta_i,\theta_j)}+F_{ij}.
\end{array}
$$
The vertical component of the Nijenhuis tensor is :
$$
\begin{array}{lll}
  \mathcal V\,N(\hat\theta_i,\hat\theta_j)&=&\left([R_{\centerdot ij}^\centerdot,Q]-
    P_i^lP_j^r[R_{\centerdot lr}^\centerdot,Q]-
    Q\Big(-P_i^l[R_{\centerdot lj}^\centerdot,Q]-
    P_j^r[R_{\centerdot ir}^\centerdot,Q]\Big)\right)\frac{\partial}{\partial Q}\\
  &=&\left[R\Big(\theta_i\wedge \theta_j-P\theta_i\wedge
    P\theta_j\Big)+ QR\Big(P\theta_i\wedge\theta_j+\theta_i\wedge
    P\theta_j\Big),Q\right]\frac{\partial}{\partial Q}\\
  &=&G(\theta_i,\theta_j)\frac{\partial}{\partial Q}.\quad\square
\end{array}
$$

\quad\\
In order to prove Theorem~1 we need to study the tensor
$G$ and we set:
\[
\left\{
  \begin{array}{l}
    G_1(\theta_i,\theta_j,P)=\theta_i\wedge \theta_j-P\theta_i\wedge
    P\theta_j\\
    G_2(\theta_i,\theta_j,P)=P\theta_i\wedge\theta_j+\theta_i\wedge
    P\theta_j.
  \end{array}
\right.
\]
An easy computation gives the following lemma.

\quad\\
{\bf Lemma 2.} {\it Let $(\theta_1,\dots,\theta_4)$ be an
  oriented orthonormal frame over an open set $\mathcal U$
  and $(I,J,K)$ be the associated basis of
  $\bigwedge^+$. Then we
  have:
  \[
  \begin{array}{l}
    I=G_1(\theta_1,\theta_2,J)=G_1(\theta_1,\theta_2,K)\\
    J=G_1(\theta_1,\theta_3,I)=G_1(\theta_1,\theta_3,K)\\
    K=G_1(\theta_1,\theta_4,I)=G_1(\theta_1,\theta_4,J)\\
    0=G_1(\theta_1,\theta_2,I)=G_1(\theta_1,\theta_3,J)=G_1(\theta_1,\theta_4,K)\\
    G_1(\theta_1,\theta_2,aI+bJ+cK)=(1-a^2)I-abJ-acK\\
    G_2(\theta_i,\theta_j,P)=PG_1(\theta_i,\theta_j,P).
  \end{array}
  \]
}

\subsection*{B)  The case  where $f$ is the  identity}
In this section we give a proof of (the well known) part~B
of Theorem~1:

\quad\\
{\bf Theorem B \cite{AHS78}.}  {\it The complex structure
  $\mathbb J_{Id}$ is integrable if and only if  $A$ is a
  homothety.}

\quad\\
The fact that $A$ is a homothety is equivalent to saying
that the selfdual Weyl tensor $W^+$ is zero. In that case
the metric is said to be {\em anti-selfdual}.

\quad\\
{\bf Proof.} In the local trivialization $\pi^{-1}(\mathcal
U)\simeq\mathcal U\times \C P^1$ of the previous section
the morphism $f=Id$ when restricted to fibers is a
holomorphic map, which only depends on the second
variable. By Proposition~1 we know that it is sufficient to
study $N(\hat\theta_i,\hat\theta_j)$. We have:
\[
\begin{array}{llll}
  F_{ij}&=&-Q_i^r[\Gamma_{r\centerdot}^\centerdot,Q]\frac{\partial}{\partial
    Q}Q_j^r\;\;\hat\theta_r
  +Q_j^r[\Gamma_{r\centerdot}^\centerdot,Q]\frac{\partial}{\partial
    Q}Q_i^l\;\;\hat\theta_l\\
  &&-Q\Big([\Gamma_{j\centerdot}^\centerdot,Q]\frac{\partial}{\partial
    Q}Q_i^l\;\;\hat\theta_l
  -[\Gamma_{i\centerdot}^\centerdot,Q]\frac{\partial}{\partial
    Q}Q_j^r\;\;\hat\theta_r\Big)\\
  &=&-Q_i^r[\Gamma_{r\centerdot}^\centerdot,Q]\hat\theta_j
  +Q_j^r[\Gamma_{r\centerdot}^\centerdot,Q]\hat\theta_i
  -Q\Big([\Gamma_{j\centerdot}^\centerdot,Q]\hat\theta_i
  -[\Gamma_{i\centerdot}^\centerdot,Q]\hat\theta_j\Big).\\
\end{array}
\]
Using
$[\Gamma_{i\centerdot}^\centerdot,Q]=[\nabla_{\theta_i}\centerdot,Q]
=\nabla_{\theta_i}
Q$ one gets:
\[
\begin{array}{lll}
  d\pi(F_{ij}) &=&-(\nabla_{Q\theta_i}
  Q)\theta_j+(\nabla_{Q\theta_j} Q)\theta_i-Q\Big((\nabla_{\theta_j}
  Q)\theta_i-(\nabla_{\theta_i}
  Q)\theta_j\Big)\\
  &=&-\nabla_{Q\theta_i}\;Q\theta_j+Q\nabla_{Q\theta_i}\;\theta_j+
  \nabla_{Q\theta_j}\;Q\theta_i-Q\nabla_{Q\theta_j}\;\theta_i\\
  &&-Q\nabla_{\theta_j}\;Q\theta_i-\nabla_{\theta_j}\theta_i+
  Q\nabla_{\theta_i}\;Q\theta_j+\nabla_{\theta_i}\theta_j\\
  &=&-E(\theta_i,\theta_j).
\end{array}
\]
The horizontal component of $N(\hat\theta_i,\hat\theta_j)$
is then zero. The vertical component is:
$$
G(\theta_i,\theta_j)=\Big[R\big(\theta_i\w
\theta_j-Q\theta_i\w Q\theta_j\big)+Q
R\big(\theta_i\w Q\theta_j+Q\theta_i\w \theta_j\big),Q\Big].
$$

\noindent But $Q$ preserves the orientation, hence:
\[
\left\{
  \begin{array}{ccc}
    \theta_i\w \theta_j-Q\theta_i\w Q\theta_j&\in&\bw^+T_mM\\
    \theta_i\w Q\theta_j+Q\theta_i\w \theta_j&\in&\bw^+T_mM.
  \end{array}
\right.
\]
Recall that the matrix of the curvature operator ${R}$ has the
following splitting:
\begin{displaymath}
  {R}=\left( \begin{array}{cc}
      A& ^tB\\
      B&C \end{array} \right)
\end{displaymath}
Since the elements of $\bw^+$ of $\bw^-$ commute \cite{AHS78}, the
component $A$ in the matrix $R$ is the only one which matters in
the computation of $G(\theta_i,\theta_j)$. By Lemma~2, one has the
equality:
$$
\big(\theta_i\w \theta_j-Q\theta_i\w
Q\theta_j\big)+Q\big(\theta_i\w Q\theta_j+Q\theta_i\w
\theta_j\big)=0, \quad\forall \theta_i,\theta_j\in T_mM.
$$
Therefore, if the matrix $A$ is a homothety the Nijenhuis
tensor of $\mathbb J_{Id}$ is zero.

Conversely, assume that $\mathbb J_{Id}$ is integrable. We
have noticed that the orthonormal frame
$(\theta_1,\dots,\theta_4)$ over $\mathcal U$ defines an
oriented orthonormal basis $(I,J,K)$ of $\bw^+$ over
$\mathcal U$.  Since $G(\theta_i,\theta_j) =0$ for all
$i,j\in\{1,\dots,4\}$, Lemma~2  implies:
\[
\begin{array}{llll}
  \textrm{at  the point } (m,I),&
  G(\theta_1,\theta_3)=&[A(J)+IA(K),I]&=0\\
  \textrm{at the point }
  (m,J),&
  G(\theta_1,\theta_2)=&[A(I)+JA(-K),J]&=0\\
  \textrm{at the point } (m,K),&
  G(\theta_1,\theta_2)=&[A(I)+KA(J),K]&=0.\\
\end{array}
\]
Since $(I,J,K)$ is an oriented orthonormal basis, it follows
from $IJ=-JI=K$ that relations of the following type hold:
$$[A(J),I]=2<A(J),K>J-2<A(J),J>~K.$$ From the previous system
we then deduce the following one:
\[
\left\{
  \begin{array}{llllr}
    <A(J),J>&=&-<IA(K),J>&=&<A(K),K>\\
    <A(J),K>&=&-<IA(K),K>&=&-<A(K),J>\\
    <A(I),I>&=&-<JA(-K),I>&=&<A(K),K>\\
    <A(I),K>&=&-<JA(-K),K>&=&-<A(K),I>\\
    <A(I),I>&=&-<KA(J),I>&=&<A(J),J>\\
    <A(I),J>&=&-<KA(J),J>&=&-<A(J),I>
  \end{array}
\right.
\]
But the matrix $A$ in the basis $(I,J,K)$ is symmetric, thus
$A$ is a homothety.  \qed 

\subsection*{C)  The case when  $f$  is constant}
\subsubsection*{Integrability theorem}
In this section we give a proof of part~C of Theorem 1.

\quad\\
{\bf Theorem C.}  {\it Let $(M,g,J_M)$ be an almost complex
  manifold such that $J_M$ is compatible with the
  orientation and the metric.  The complex structure
  $\mathbb J_\infty$ is integrable if and only if:
\begin{enumerate}
 \item[{\rm i)}] $J_M$ is integrable;

  \item[{\rm ii)}] the kernel of $A$ contains the subspace
  $J_M^\perp\subset \bigwedge^+$ orthogonal to the line
  generated by $J_M$ $($i.e. $J_M^\perp\subset \ker(A))$.
\end{enumerate}
}
\quad\\
Notice that the integrability condition is not conformal on $g$.
Moreover, when $\mathbb J_\infty$ is integrable, it gives to the
twistor projection $\pi : (Z,\mathbb J_\infty)\lra (M,J_M)$ the
structure of a holomorphic $\C P^1$-bundle.

For a complex manifold $(M,g,J_M)$ we have a decomposition
$\C\otimes TM=T^{1,0}\oplus T^{0,1}$ into $\pm i$
eigenspaces of $J_M$. We then obtain:
\[
\left\{
  \begin{array}{l}
    \C\otimes \bigwedge^+=\C
    J_M\oplus^\perp(\bigwedge^{2,0}\oplus\overline{\bigwedge^{2,0}})\\
    \C\otimes\bigwedge^-=\{\psi\in\bigwedge^{1,1} \, \mid \,
    <\psi, J_M>=0\}
  \end{array}
  \textrm{ where } \left\{
    \begin{array}{l}
      \bigwedge^{2,0}=T^{1,0}\wedge T^{1,0}\\
      \bigwedge^{1,1}=T^{1,0}\wedge T^{0,1}
    \end{array}
  \right.  \right.
\]
\quad\\
Condition ii) says that
$(\bigwedge^{2,0}\oplus\overline{\bigwedge^{2,0}})\subset
\ker(A)$.  For a K\"ahlerian manifold the curvature $R$ may
be viewed as a symmetric endomorphism of $\bigwedge ^{1,1}$,
so in some orthonormal basis compatible with these
decompositions we have
$A=\left[\begin{array}{lll}\frac{s}{4}&0&0\\
    0&0&0\\
    0&0&0
  \end{array}
\right]$ and $W^+=\left[\begin{array}{ccc}\frac{s}{6}&0&0\\
    0&-\frac{s}{12}&0\\
    0&0&-\frac{s}{12}
  \end{array}
\right]$. We then have the following result:

\smallskip

\quad\\
{\bf Proposition 3.}  {\it For any K\"ahlerian manifold
  $(M,g,J_M)$ the complex structure $\mathbb J_\infty$ on
  $(Z,\tilde g)$ is integrable. Furthermore, if $(M,g,J_M)$
  is K\"ahler and the scalar curvature of $g$ is never zero,
  then $\mathbb J_\infty$ and $\mathbb J_{-\infty}$ (the
  compatible complex structure on $(Z,\tilde g)$ associated
  to $-J_M$) are the only compatible complex structures on
  $(Z,\tilde g)$}.

\quad\\
The proof will show that the result is locally true. In other
terms, for a K\"ahlerian manifold whose scalar curvature is non
zero there are, even locally, only two compatible complex
structures on its twistor space.

\quad\\
{\bf Proof.} The first part being a consequence
 of Theorem~C,
we only need to prove the second part of the proposition. Let
$\mathbb J_f$ be a compatible complex structure on $(Z,\tilde g)$
and assume that the scalar curvature of $(M,g,J_M)$ is never zero.
One can build an orthonormal basis $(I,J,K)$ of $\bigwedge^+$ over
an open set $\mathcal U$ as follows. Setting $I=J_M$, pick any
unitary vector $J$ orthonormal to $I$ and define $K=IJ$. For any
$m\in \mathcal U$, there exists $(a,b,c)\in\Sp^2$ such that
$f(m,J)=aI+bJ+cK$. But, as $(M,g,J_M)$ is K\"ahler, in this basis
we have
$A=\left[\begin{array}{lll}\frac{s}{4}&0&0\\
    0&0&0\\
    0&0&0
  \end{array}
\right]$.  Let $\theta_1$ be a unitary vector field defined
over $\mathcal U$; set $\theta_2=I\theta_1$. As $\mathbb
J_f$ is integrable,  $G(\theta_1,\theta_2)$ is identically zero
on $\mathcal U$. In particular, at the point
$(m,J)$ we obtain:
\[
\begin{array}{lll}
  G(\theta_1,\theta_2)&=&0\\
  &=&[A\Big((1-a^2)I-abJ-acK\Big)+JA(cJ-bK),J]=0\\
  &=&[(1-a^2)\frac{s}{4}I,J]=(1-a^2)\frac{s}{2}K.
\end{array}
\]
Therefore $a=\pm 1$, that is $f(m,J)=\pm I$ for all $J$
orthonormal to $I$. Since $f$ must be holomorphic in
the  fibers we get that  $f$ is constant, equal to $I$ or
$-I$. \qed

\quad\\
{\bf Proof of Theorem~C.} By Proposition~1, it is sufficient to
check that $N(\hat\theta_i,\hat\theta_j) =0$. As $f$ is constant
on fibers we always have $F_{ij}=0$. Therefore:  $\mathbb
J_\infty$ integrable $\iff$
$E(\theta_i,\theta_j)=G(\theta_i,\theta_j)=0$ $\iff$ \{$J_M$
integrable and $G(\theta_i,\theta_j)=0$\}. But for all
$\theta_i,\theta_j\in TM$ we have $\left\{\begin{array}{l}
    \theta_i\wedge \theta_j-J_M\theta_i\wedge J_M\theta_j
    \in J_M^\perp\\
    J_M\theta_i\wedge \theta_j+\theta_i\wedge J_M\theta_i\in
    J_M^\perp\end{array}\right.$. Consequently,  if $J_M^\perp\subset
\ker(A)$ we obtain $G(\theta_i,\theta_j)=0$ for all
$\theta_i,\theta_j\in TM$.

Conversely, suppose that $\mathbb J_\infty$ is integrable. Set $J_0=J_M$. Locally
over an open set $\mathcal U$ one can complete $\{J_0\}$ to get an
oriented orthonormal basis $(I_0,J_0,K_0)$ of $\bw^+$. Let
$\theta_1$ be a unitary vector field defined over $ \mathcal U$;
set $\theta_2=I_0\theta_1$. If $G=0$, then, for all $m\in \mathcal
U$ and $Q\in \pi^{-1}(m)$, Lemma~2 implies that at the point
$(m,Q)$:
\[
\begin{array}{lcl}
  G(\theta_1,\theta_2)=&[A(I_0)+QA(-K_0),Q]&=0.
\end{array}
\]
In particular, for $Q=A(K_0)$, we have $[A(I_0),A(K_0)]=0$
and it follows that $A(K_0)=cA(I_0)$ for some constant
$c$. The former equation yields:
$$
\forall Q\in
\pi^{-1}(m), \quad 0=[A(I_0)+QA(-K_0),Q]=(Id-cQ)\,[A(I_0),Q]
\Longrightarrow A(I_0)=0.
$$
Therefore $J_0^\perp=\Vect(I_0,K_0)\subset \ker \,A$. \qed

\quad\\
Recall that we have a characterization of an integrable
almost complex structure $J_M$ on $M$ in terms of the twistor
space and one of the K\"ahlerian complex structures.

\quad\\
{\bf Proposition (see, for example, \cite{Sal85, BdB88}).}
{\it Let  $J_M$ be a Hermitian almost complex structure  on
  $(M,g)$.
  Then:
\begin{itemize}
\item $J_M$ is integrable if and only if the associated
  section of the twistor space, $s: (M,J_M) \lra (Z,\mathbb
  J_{Id})$, is almost holomorphic, that is: the differential
  $ds$ satisfies $ds\circ J_M=\mathbb J_{Id}\circ ds$;
\item $J_M$ is K\"ahlerian if and only if $s$ is an
  horizontal section, that is to say: the tangent space of
  the submanifold $s(M)\subset Z$ is included in the
  horizontal distribution.
\end{itemize}
  }

  \quad\\
  It is well known that the existence of a K\"ahlerian metric
  on a compact complex surface $(M,J_M)$ is equivalent for
  the first Betti number $b_1$ to be even \cite{Mio74,
    Siu83, Lam99}. Theorem~C gives a new characterization of
  compact K\"ahlerian manifolds in terms of compatible
  complex structures on the associated twistor spaces.

  \quad\\
  {\bf Proposition 4.} {\it A compact almost Hermitian
    $4$-dimensional manifold $(M,g,J_M)$ is K\"ahlerian if
    and only if $\mathbb J_\infty$ is integrable.}

\quad\\
In section~E we will deduce from  that  proposition a
characterisation of compact scalar-flat K\"ahler manifolds in
terms of compatible complex structures on $(Z,\tilde g)$
(cf.~Proposition~8).

\quad\\
{\bf Proof.} Let $\theta$ be the Lee form of $(M,g,J_M)$ defined
by $dJ_M=-2\theta\wedge J_M$, where $J\in\bigwedge^+$ is viewed as
a 2-form. Denote by $\kappa$ the conformal scalar curvature, which
is related to the scalar curvature $s$ by $\kappa=s+6(\delta
\theta-\vert\theta\vert^2)$. The condition $J_M^\perp\subset \ker
A$ is equivalent to the following: the selfdual Weyl tensor $W^+$
is degenerate (meaning that, in every point,  two of the
eigenvalues coindident) and the scalar curvature of $(M,g)$ is
equal to the conformal scalar curvature \cite{AG08}. This is also
equivalent to $\delta\theta=\vert\theta\vert^2$. Integrating this
expression over $M$ gives $\theta=0$ by the Brochner-Grenn
theorem. But $(M,g,J_M)$ is K\"ahler if and only if $\theta$
vanishes identically. \qed

\quad\\
{\bf Corollary 1.} {\it Assume that a compact $4$-dimensional
manifold $(M,g)$ admits two almost complex structures $J_1\neq \pm
J_2$ compatible with the metric and the orientation. Then the
associated compatible almost complex structures $\mathbb
J_{\infty1}, \mathbb J_{\infty2}$ on $(Z,\tilde g)$ are integrable
if and only if $\{J_1, J_2\}$ spans a hyperk\"ahler structure on
$(M,g)$.}

\quad\\
{\bf Proof.} By Proposition~4, $\mathbb J_{\infty1}$ and $\mathbb
J_{\infty2}$ are integrable if and only if $J_1$ and $J_2$ are
K\"ahler. As $J_1\neq\pm J_2$, then $J_1$ is different from $\pm
J_2$ everywhere. The holonomy of $g$ reduces to $U(2)$ by $J_1$
and further to $SU(2)$ by $J_2$. This says that $g$ is
hyperk\"ahler. \qed


\subsubsection*{Study of the manifold $(Z,\mathbb J_\infty)$}

Any scalar-flat K\"ahler manifolds $(M,g,J_M)$ is automatically
anti-selfdual \cite{Gau81}. For such a manifold we can put two
natural complex structures on its twistor space: $\mathbb J_{Id}$
and $\mathbb J_\infty$. The next proposition shows that these
complex structures are never deformation of each other.

\quad\\
{\bf Proposition 5.} {\it If $(M,g,J_M)$ is a scalar-flat
  K\"ahler manifold, the complex structure $\mathbb J_\infty$
  on $Z$ is never a deformation of the complex structure
  $\mathbb J_{Id}$.}

\quad\\
{\bf Proof.} It is sufficient to show that $(Z,\mathbb
J_{Id})$ and $(Z,\mathbb J_\infty)$ do not have the same
Chern classes. Let $h$ be the generator of the second
cohomology group $H^2(\C P^1,\Z)\simeq\Z$. By Leray-Hirsch
theorem's \cite{BT82} the  cohomology ring of $Z$ is a
$H^\star(M,\R)$-module generated by $h$ with
relation $4h^2=3\tau+2\chi$, where $\tau$ and $\chi$ are the
signature and the Euler characteristic of $M$. Denote by
$c_1(J_M)$ the first Chern class of the manifold $(M,J_M)$.
Under this notation we have :
\[\begin{array}{ll}
  c(\mathbb J_{Id})&=1+4h+3\tau+3\chi+2h\chi \quad\cite{Hit81}\\
  c(\mathbb J_\infty)&=(1+2h)(1+c_1(J_M)+\chi)\\
  &=1+2h+c_1(J_M)+2hc_1(J_M)+\chi+2h\chi.
\end{array}
\]
If the complex structures were deformations of each other,
they would have the same Chern numbers: $c_1(\mathbb
J_{Id})^3=16(3\tau+2\chi)h=c(\mathbb
J_\infty)^3=8(3\tau+2\chi)h$.  This forces $3\tau+2\chi =0$.
Let $\mu_g$ be the volume form on $M$ associated to the
metric $g$; by the Gauss-Bonnet formula \cite{AW43},
\cite{Hir66}:
$$
3\tau+2\chi=\frac{1}{4\pi^2}\int_M 2\Vert
W^+\Vert+\frac{1}{24}s^2-2\Vert
B\Vert^2\mu_g=-\frac{1}{2\pi^2}\int_M\Vert B\Vert^2\mu_g.
$$
Thus, $3\tau+2\chi=0$ implies $B=0$. As the scalar curvature of
$(M,g)$ is supposed to be zero, the manifold $(M,g,J_M)$ would be
Ricci-flat, hence $c_1(J_M)=0$. Therefore the first Chern classes
of $(Z,\mathbb J_{Id})$ and of $(Z,\mathbb J_\infty)$ are
different and these two manifolds are never
deformations of each other. \qed\\

When $(M,g,J_M)$ is a complex spin manifold, Hitchin has
shown that there exists a holomorphic line bundle $L\lra M$
such that $L\otimes L=K_M$ is the canonical line bundle
\cite{Hit74b}. Then, the twistor space $Z$ can be
identified, in a $\mathcal C^\infty$-way, to the
projectivization bundle $\PP(L\oplus L^\star)$
\cite{Sal83}. By this construction we see that the manifold
$Z\simeq \PP(L\oplus L^\star)$ admits a natural complex
structure denoted by $\mathbb I$. When $(M,g,J_M)$ is not
spin, but only complex, the bundle $L\oplus L^\star$ exists
only locally. Nevertheless, the projectivization
$\PP(L\oplus L^\star)$ still exists globally, due to the
fact that the transition functions on $L\oplus L^\star$ are
well defined holomorphic maps up to sign. In general
$\mathbb I$ is not a compatible complex structure on
$(Z,\tilde g)$.

Now, if $(M,g,J_M)$ satisfies the conditions of Theorem~C, we can
put another complex structure on its twistor space, namely
$\mathbb J_\infty$. The question is then to determine the
relationship between the manifolds $(Z,\mathbb I)$ and $(Z,\mathbb
J_\infty)$. In that direction we have the following result.

\quad\\
{\bf Proposition 6.} {\it Let $(M,g,J_M)$ be a manifold
  satisfying conditions of Theorem~C (i.e.~$\mathbb
  J_\infty$ integrable). The complex structures $\mathbb I$
  and $\mathbb J_\infty$ on $Z$ are deformations of each
  other: there exists on $Z$ a path of
  integrable complex structures $\mathbb J_t$, $t\in[0,1]$,
  connecting $\mathbb I$ to $\mathbb J_\infty$.}

\quad\\
By combining this result 
and \cite[Theorem~4.1]{Tsa87} we obtain
 another proof of Proposition~5.

\quad\\
{\bf Proof.} In an appropriate local trivialization of the bundle
$Z\lra M$, the almost complex structure $\mathbb I$ on $\mathcal
U\times \Sp^2$ can be identified with the product structure
$J_M\times J_{\C P^1}$. Let
$(\theta_1,\theta_2,\theta_3,\theta_4)$ be an oriented orthonormal
frame defined over $\mathcal U$ providing this trivialization. Set
$\hat\theta_{i,t}=\theta_i-t[\Gamma_{i\centerdot}^\centerdot,Q]\frac{\partial}
{\partial Q}$ for $t\in[0,1]$. The subspace
$H_t=\Vect(\hat\theta_{1,t},\dots,\hat\theta_{4,t})$ is in direct
sum with the vertical distribution $V_z$ and can be glued into a
global distribution over $Z$. Define the almost complex structure
$\mathbb J_t$ on $\pi^{-1}{\mathcal
  U}$ as follows: endow $V_z$ with the complex structure
of the fibers (complex projective lines) and  pull back on
$H_t\simeq T_m M$ the complex structure $J_M$. Then, $\mathbb
J_t$ is a path of almost complex structures from $\mathbb I$
to $\mathbb J_\infty$. The integrability of $\mathbb J_t$ is
shown in the same way as that of $\mathbb J_\infty$. \qed

\subsection*{D) The case where $f$ is a homothety}
\subsubsection*{Integrability theorem}
In this section we prove part~D of Theorem~1.

\quad\\
{\bf Theorem D.} {\it Let $(M,g,J_M)$ be a K\"ahlerian
  manifold. For all complex $\lambda\notin\{0,1\}$ the
  almost complex structure $\mathbb J_{\lambda Id}$ is
  integrable if and only if $(M,g,J_M)$ is scalar-flat
  K\"ahler $($i.e.~$A=0)$.}

\quad\\
The condition $A=0$ is equivalent to saying that the metric $g$ is
Hermitian scalar-flat and anti-selfdual. These metrics are called
{\it optimal} by LeBrun because they are absolute minimizers of
the functional $\mathcal K(g)=\int_M\vert R\vert^2 d \textrm{vol}$
\cite{Leb08}. Let $(M,g,J_M)$ be a compact scalar-flat K\"ahler
manifold and $c_1(M)$ be the real first Chern class of $(M,J_M)$.
Two possibilities may occur \cite{Laf82}. Either $c_1(M)=0$ and
($M,g,J_M)$ is then finitely covered by a hyperk\"ahler surface,
i.e.~a flat torus or a $K3$-surface with Ricci-flat K\"ahler
metric \cite{Boy86}, \cite{Pon91a}. Or $c_1(M)<0$, in which case
$(M,g)$ is a ruled surface \cite{KLP97}, i.e.~$(M,g)$ is obtained
by blowing up $m$ points on a $\C P^1$-bundle over a Riemann
surface of genus $\gamma$. The condition $c_1(M)<0$ gives a lower
bound on the number of points $m$ to be blown up: namely $m\geq 9$
when $\gamma=0$, $m\geq 1$ when $\gamma=1$ and there is no
restriction for $\gamma>1$. Conversely we have:

\quad\\
{\bf Theorem \cite{KLP97}.} {\it A ruled surface $M$ has a
  blow-up $\tilde M$ which is a scalar-flat K\"ahler
  surface. Moreover, any further blow up of $\tilde M$ admits
  a scalar-flat K\"ahler metric.}

\quad\\
For simply connected manifold we have the following result:

\quad\\
{\bf Theorem \cite{KLP97}, \cite{Leb08}.} {\it Let $M$ be a
  smooth compact simply connected $4$-manifold. If $M$
  admits a scalar-flat K\"ahler structure, then $M$ is
  diffeomorphic to a $K3$-surface or to the connected sum
  $\C P^2\sharp k\overline{\C P^2}$ for some
  $k\geq10$. Conversely, if $M$ is a $K3$-surface or is
  diffeomorphic to $\C P^2\sharp k\overline{\C P^2}$ for
  some $k\geq14$, then it admits a scalar-flat K\"ahler
  metric.}

\quad\\
{\bf Proof of Theorem~D.} By Propositions~1~\&~2, if $A=0$ it is
enough to show that $E(\theta_i,\theta_j)+F_{ij}=0$ to get the
integrability of $\mathbb J_{\lambda Id}$.  Let $z\in \pi^{-1}(m)$
be a point of $Z$ over $m\in M$. Let $\theta_1,\theta_2$ be two
unitary vector fields, defined on an open set $\mathcal U$ of $M$,
such that $\nabla\vert_m\theta_1=\nabla\vert_m\theta_2=0$ and
$\theta_2\in (\theta_1,J_M\theta_1)^\perp$. As $J_M$ is parallel,
the vector fields $\theta_3=J_M\theta_1$ and
$\theta_4=J_M\theta_2$ complete the family $(\theta_1,\theta_2)$
to give an orthonormal basis such that $\nabla\vert_m\theta_i=0$
for all $i\in\{1,\dots,4\}$. Hence $F_{ij}\vert_m=0$, since
$\Gamma_{ij}^k(m)=0$ for all $i,j,k\in\{1,\dots,4\}$. Moreover,
this frame gives  an oriented orthonormal basis $(I,J,K)$ of
$\bigwedge^+$, and therefore a local trivialization of $Z$ over
the open set $\mathcal U$, where $\infty$ coincides with $J$. It
follows that the restriction of $f$ to the fibers does not depend
on the second variable:
\[
\xymatrix @R=0,2pc @C=0,5cm{ \pi^{-1}(\mathcal U)\simeq
  \mathcal U\times\Sp^2&\ar[rr]^{f}&&&
  \pi^{-1}(\mathcal U)\simeq \mathcal U\times\Sp^2\\
  (x,Q)\ar[ddddr]_{\pi}&\ar[rr]&&& \Big(x,f(Q)\Big)\ar[ddddl]_{\pi}\\
  \\
  \\
  \\
  && \mathcal U&&\\
}
\]
Thus $E(\theta_i,\theta_j)\vert_m=0$ and $\mathbb J_{\lambda
  Id}$ is integrable.

Conversely, assume that $\mathbb J_{\lambda Id}$ is integrable.
From Proposition~3 one deduces that the scalar curvature must be
zero, hence $A=0$. \qed


\subsubsection*{Study of the manifold $(Z,\mathbb J_{\lambda
    Id})$}
We know that the almost complex structure $\mathbb J_{Id}$
on $Z$ is integrable if and only if the metric $g$ is
anti-selfdual. In that case the twistor space $Z$ is a
complex $3$-manifold. The fibers of the projection $\pi:
Z\lra M$ are rational curves with normal bundles isomorphic
to $\mathcal O(1)\oplus \mathcal O(1)$.  On each fiber the
antipodal map $\sigma$ is an antiholomorphic free
involution. Observe that this construction only depends on
the conformal class of the metric $g$. The converse holds:
an arbitrary conformal anti-selfdual $4$-manifold can be
constructed from a complex $3$-manifold $Z$ as soon as $Z$
admits a free antiholomorphic involution $\sigma$ and a
foliation by $\sigma$-invariant rational curves, each of
which having $\mathcal O(1)\oplus \mathcal O(1)$ as normal
bundle \cite{Pen76}, \cite{AHS78}.  This is the Penrose
correspondence.

Thus if $(M,g,J_M)$ is a scalar-flat K\"ahler surface we have
that, for $\lambda\in\C-\{0,1\}$, the complex $3$-manifold
$(Z,\mathbb J_{\lambda Id})$ is the twistor space of
$(M,g_\lambda)$ for some anti-selfdual metric $g_\lambda$.

At least two cases may occur. Firstly: all the $(Z,\mathbb
J_{\lambda Id})$ are biholomorphic to $\mathbb J_{Id}$, thus there
exists a $1$-dimensional family of biholomorphism of $(Z,\mathbb
J_{Id})$. We will see in section~G that this is the case for any
bi-elliptic surface (quotient of a flat torus). Secondly: none of
the $(Z,\mathbb J_{\lambda Id})$ is biholomorphic, thus we have a
$1$-dimensional family of non conformal anti-selfdual metrics on
$M$. For example, if one blows-up at least $14$ points in $\C
P^2$, equipped with the Fubini-Study metric $g$, one gets $(\C
P^2\sharp k\overline{\C P^2},g)$ for some $k\geq 14$. This
manifold admits a scalar-flat K\"ahler metric \cite{Leb08} but
doesn't have any non trivial conformal application, thus its
twistor space doesn't have any biholomorphism. Therefore the
complex structure $(Z,\mathbb J_{\lambda Id})$ defines a
$1$-dimensional family of anti-selfdual metrics:

\quad\\
{\bf Proposition 7.} {\it There exist a $1$-dimensional
  family of non conformal scalar-flat K\"ahler metrics on $(\C
  P^2\sharp k\overline{\C P^2},J_M)$ for every $k\geq 14$.}

\quad\\
{\bf Proof.} There exists a $1$-dimensional family of non
conformal anti-selfdual metrics $g_\lambda$ on $(\C P^2\sharp
k\overline{\C P^2},J_M)$. But $g_\lambda$ is Hermitian, thus in
the conformal class of $g_\lambda$ there exists a scalar-flat
K\"ahler metric \cite{Boy88}. \qed

\subsection*{E)  Metric properties on $M$ in terms of
  compatible complex structures on $(Z,\tilde g)$}
We can use the almost complex structures $\mathbb J_f$ to
characterize some properties of the metric $g$ on $M$. Indeed, by
(the well known) Theorem~B we have that $g$ is anti-selfdual if
and only if $\mathbb J_{Id}$ is integrable. We  showed that a
compact almost Hermitian manifold $(M,g,J_M)$ is K\"ahler if and
only if $\mathbb J_\infty$ is integrable; furthermore the
integrability of $\mathbb J_{Id}$ and $\mathbb J_\infty$ is
equivalent to $(M,g,J_M)$ scalar-flat K\"ahler
(cf.~Proposition~8).

When limiting to the case where $(M,g)$ is anti-selfdual, we
can give a characterization of metrics which are scalar-flat
in terms of compatible complex structures on $(Z,\tilde g)$.
According to the terminology of LeBrun this is a
characterization of optimal metrics \cite{Leb08}.

\quad\\
{\bf Theorem E.} {\it Let $(M,g)$ be an anti-selfdual
  Riemannian manifold. The following are equivalent:
\begin{itemize}
\item the scalar curvature of $g$ is flat;
\item every $ m\in M$ has an open neighborhood $\mathcal U$
  such that $Z$ admits, over $\mathcal U$, an integrable
  compatible complex structure $\mathbb J_f$ for at least
  one (and then infinitely many) morphism(s) $f\neq Id$.
\end{itemize}
}
\quad\\
 In other words, if $(M,g)$ is an anti-selfdual metric with
non zero scalar curvature then, even locally on $Z$, the only
integrable almost complex structure among the $\mathbb J_f$'s is
$\mathbb J_{Id}$. This result should be compared to the following
result of Salamon:

\quad\\
{\bf Proposition \cite{Sal91} (see also \cite{Pon97}).} {\it
  A metric $g$ on $M$ is anti-selfdual if, and only if,
  locally around each point $m\in M$ there exist infinitely
  many compatible complex structures on $(M,g)$.}

\quad\\
In a similar direction, Pontecorvo gives a conformal
characterization of scalar-flat K\"ahler manifolds among
anti-selfdual Hermitian manifolds.  Indeed, let $(M,g,J_M)$ be an
anti-selfdual complex Hermitian manifold. The complex structure
$J_M$ on $M$ defines a section $s :Z\lra M$ \cite{BdB88}, whose
image will be noted $\Sigma=s(M)$. Similarly, the hypersurface
$\overline{\Sigma}=\sigma(\Sigma)$ of $Z$ corresponds to the
conjugate complex structure $-J_M$. Let $X$ be the divisor
$\Sigma+\overline{\Sigma}$ in $Z$ and consider the holomorphic
line bundle $[X]$. Denote by $K_Z$ be the canonical line bundle of
$(Z,\mathbb J_{Id})$.

\quad\\
{\bf Proposition \cite{Pon92a}.}  {\it Let $(M,g,J_M)$ be a
  Hermitian anti-selfdual manifold. The line bundles $[X]$
  and $-\frac{1}{2}K_Z$ are isomorphic if and only if $g$ is
  conformal to a scalar-flat K\"ahler metric.}

\quad\\
Notice that Theorem~1 and Proposition~3\&4 give a non conformal
characterization of compact scalar-flat K\"ahler manifolds.

\quad\\
{\bf Proposition 8.}  {\it Let $(M,g,J_M)$ be a compact
  almost Hermitian manifold. The following are equivalent:
\begin{itemize}
\item the metric $g$ is scalar-flat K\"ahler;
 \item the compatible
complex structures $\mathbb
  J_{Id}$ and $\mathbb J_{\infty}$ on $(Z,\tilde g)$ are
  integrable;
\item the compatible complex structures $\mathbb
  J_{\lambda Id}$ and $\mathbb J_{\infty}$ on $(Z,\tilde g)$
  are integrable.
\end{itemize}
}

\quad\\
{\bf Proof.} A K\"ahlerian manifold $(M,g,J_M)$ is scalar-flat if
and only if  $g$ is anti-selfdual \cite{Gau81}. Then, it follows
from Proposition~3\&4 and Theorem~1 that: \{$\mathbb J_{\infty}$
and $\mathbb J_{\lambda Id}$ are integrable\} $\iff$
 \{$g$ is scalar-flat K\"ahler\}
 $\iff$ \{$(M,g,J_M)$ is anti-selfdual K\"ahler\} $\iff$
\{$\mathbb J_\infty$ and $\mathbb J_{Id}$ are integrable\}. \qed

\quad\\
{\bf Proof of Theorem~E.} If $(M,g)$ is a scalar-flat
anti-selfdual metric its twistor space is complex and $(M,g)$
admits, locally, at least one complex structure $J_M$
\cite{Sal91}. Then Theorem~C ensures that the locally defined
almost complex structure $\mathbb J_\infty$ on $Z$ is integrable.
Actually, as soon as $(M,g)$ is scalar-flat there are, locally,
infinitely many integrable complex structures $J_M$ on $M$, and so
infinitely many integrable complex structures $\mathbb J_\infty$
on $Z$.

Conversely, let $(M,g)$ be a manifold with an anti-selfdual
metric $g$ having non zero scalar curvature. Let $f:Z\lra Z$
be a morphism such that $\mathbb J_f$ is integrable over an
open set $\mathcal U$. Let $(m,Q)$ be a point in
$\pi^{-1}(\mathcal U)$ and set $f(m,Q)=P$. According to our
notation, if $\mathcal U$ is small enough we can build an
orthonormal basis $(\theta_1,\dots,\theta_4)$ of vector
fields on $M$ such that $P=J=
\theta_1\wedge\theta_3-\theta_2\wedge\theta_4$. Then there
exists $(a,b,c)\in\Sp^2$ such that $Q=aI+bJ+cK$.

As $\mathbb J_f$ is integrable, $G(\theta_1,\theta_2)$
vanishes everywhere. In particular, at the point $(m,Q)$ one
obtains:
\[
\begin{array}{lll}
  G(\theta_1,\theta_2)&=&0\\
  &=&\frac{s}{12}[I-QK,Q]\\
  &=&\frac{2s}{12}\Big(acI-c(1-b)J+\big(b(1-b)-a^2\big)K\Big)
\end{array}
\Longrightarrow \left\{
  \begin{array}{l}
    ac=0\\
    c(1-b)=0\\
    b=a^2+b^2
  \end{array}
\right.
\]
Therefore we  have $Q=J=P$ for every $(m,Q)\in
\pi^{-1}(\mathcal U)$,  that is to say $f=Id$. \qed

\subsection*{F)  Compatible complex structure on locally
  conformally K\"ahler manifolds}
The aim of this section is to give a local description of the set
$\mathcal I$ of integrable compatible complex structures on the
twistor space $(Z,\tilde g)$ of a compact locally conformally
K\"ahler (abbreviated in {\it l.c.k.}) manifold $(M,g,J_M)$. This
condition is equivalent to $W^+$ being degenerate, which means
that at each point of $M$ at least two eigenvalues of $W^+$
coincide.

We start by recalling the main results about the {l.c.k.}
manifolds.

A result by Tricerri, generalizing the analogous result in the
K\"ahler case, shows that it is enough to understand minimal
complex surfaces.

\quad\\
{\bf Proposition \cite{Tri82}.} {\it A complex manifolds
  $(M,g,J_M)$ is l.c.k if and only if the blow-up of
  $M$ at a point is l.c.k.}

\quad\\
When the first Betti number $b_1$ is even, a {l.c.k.} manifold is
globally K\"ahler.

\quad\\
{\bf Proposition \cite{Vai80}.} {\it Every l.c.k.
  metric on a compact surface $(M,J_M)$ with even first
  Betti number is globally conformal K\"ahler.}

\quad\\
When the first Betti number is odd and the Euler
characteristic is zero, we have a classification due to
Belgun, Gauduchon-Ornea, Tricerri, Vaisman.

\quad\\
{\bf Proposition \cite{Bel00}.} {\it The complete list of
  compact minimal l.c.k. surfaces with odd first Betti
  number and zero Euler characteristic is:
\begin{enumerate}
 \item[{\rm i)}] the properly elliptic surfaces (i.e.~surfaces with
  $\Kod(M)=1$ and $b_1$ odd);

  \item[{\rm ii)}] the Kodaira surfaces (i.e.~surfaces with
  $\Kod(M)=0$ and $b_1$ odd);

  \item[{\rm iii)}] the Hopf surfaces;

  \item[{\rm iv)}] the Inoue-Bombieri surfaces different from
  $S^-_{n,u}$ with $u\notin \R$ {\em \cite{Tri82}}.
  \end{enumerate}}
\quad\\
When the first Betti number is odd and the Euler characteristic is
non zero, the only other possible case is that of surfaces of
class $VII$ with $0<\chi=b_2$ \cite{BHPVdV04}, for which there is
(yet) no classification. (For some existence results see
\cite{FP05}.)

Let $\mathbb J$ be a compatible almost complex structure on
$(Z,\tilde g)$. We say that $\mathbb J$ is semi-integrable if the
vertical component of the Nijenhuis tensor is zero. Denote by
$\mathcal I_\frac{1}{2}$ (resp.~$\mathcal I$) the set of
semi-integrable (resp.~integrable) compatible complex structures
on $(Z,\tilde g)$. Propositions~1 and~2 give a necessary and
sufficient condition for $\mathbb J$ to be semi-integrable, or
integrable. The set $\mathcal I$ on a l.c.k. manifold $(M,g,J)$
depends on the spectrum of the operator $A=W^++\frac{s}{12}$.  Let
$\kappa$ be the conformal curvature defined in proposition 4. Then
on an adapted basis we have :
\[
A=W^++\frac{s}{12}Id= \left[\begin{array}{ccc}
\frac{2\kappa}{12}&0&0\\
0&\frac{-\kappa}{12}&0\\
0&0&\frac{-\kappa}{12}
\end{array}
\right] +\left[\begin{array}{ccc}
\frac{s}{12}&0&0\\
0&\frac{s}{12}&0\\
0&0&\frac{s}{12}
\end{array}
\right]=\left[\begin{array}{ccc}
x&0&0\\
0&y&0\\
0&0&y
\end{array}
\right].
\]
Moreover $J_M$ is actually an eigenvector of $W^+$ for the simple
eigenvalue $\frac{\kappa}{6}$.

\quad\\
{\bf Theorem 2.} {\it Let $(M,g,J_M)$ be a compact surface l.c.k.,
 if we don't have $x=y=0$ we note
$\frac{x}{y}\in\R\cup\{\infty\}$. On an  open set $\mathcal U$ of
$M$  :
\begin{itemize}
\item[{\rm A)}] We have $x=y=0$ if, and only if, on $\mathcal U$
one of the following equivalent conditions hold:

\begin{itemize}
\item[{\rm i)}] $(M,g,J_M)$ is scalar-flat K\"ahler.

\item[{\rm ii)}]  $g$ anti-selfdual scalar-flat.

\item[{\rm iii)}]  The compatible complex structures  $\mathbb
J_{Id}$, $\mathbb J_\infty$ and $\mathbb J_{\lambda Id}$ are
integrable.

\item[{\rm iv)}] The cardinal of $\mathcal I$ is infiny.
\end{itemize}
This is the case globally if, and only if,  $(M,g,J_M)$ is a flat
torus (or a quotient), a $K3$-surface with a Calabi-Yau metric (or
a quotient), a $\C P^1$-bundle over a Riemann surface
$\Sigma_\gamma$ of genus $\gamma>1$ with the conformally flat
K\"ahler metric which  locally is a product of the (+1)-curvature
metric on $\C P^1$ and (-1)-curvature metric on $\Sigma_\gamma$
\cite{Boy88}, \cite{Pon92b}.

\item[{\rm B)}] We have $\frac{x}{y}=\infty$ if, and only if,  on
$\mathcal U$ one of the following equivalent conditions hold:

\begin{itemize}
\item[{\rm i)}] $(M,g,J_M)$ is K\"ahler with $s\neq 0$.

\item[{\rm ii)}] $\mathcal I=\mathcal I_{\frac{1}{2}}=\{ \mathbb
J_{-\infty},\mathbb J_\infty\}$.
\end{itemize}
This is the case globally on $M$ if $(M,g,J_M)$ is
K\"ahler-Einstein not Ricci-flat (that is a Fano manifolds or a
manifold where the canonical line bundle is ample).

 \item[{\rm C)}] We have $\vert\frac{x}{y}\vert\leq1$ if, and only if,  on $\mathcal U$:
$\mathcal I_\frac{1}{2}=\{\mathbb J_{e^{\pm i\theta}Id}\}$ where
$\cos\theta=\frac{x}{y}$.

 \item[ {\rm D)}] We have $\infty\neq\vert\frac{x}{y}\vert \geq1$
if, and only if,  on $\mathcal U$:
 $\mathcal I_\frac{1}{2}=\{\mathbb J_{u_1 Id},\mathbb J_{u_2
Id}\}$ where $u_1=\frac{1+\sin\theta}{\cos\theta}$,
$u_2=\frac{1-\sin\theta}{\cos\theta}$ and
$\cos\theta=(\frac{x}{y})^{-1}$.
\end{itemize}
}
\quad\\
{\bf Remark.} We have $\frac{x}{y}=1$ if, and only if,
$(M,g,J_M)$ is anti-selfdual with $s\neq0$. If it is the case
globally then $(M,J_M)$ must be in class VII \cite{Boy88}. We can
find some global example of manifolds $(M,g,J_M)$ with arbitrary
$\frac{x}{y}$ in \cite{AM99}.

 \quad\\
 {\bf Proof of~A.}  The multiplicity of the eigenvalue $0$ of $A$ is equal to
 $3$ $\iff$ $\kappa=s=0$ $\iff$ $(M,J_M,g)$ scalar-flat K\"ahler $\iff$
  $(M,J_M,g)$ anti-selfdual scalar-flat
 \cite{Boy88} $\iff$ $\mathbb
J_{Id}$, $\mathbb J_\infty$ and $\mathbb J_{\lambda Id}$
integrable by proposition~8.  The equivalence with
 condition iv) will be a consequence of (the rest of the proof of) the
theorem.   \qed

 \quad\\
 {\bf Proof of B.} The multiplicity of the eigenvalue $0$ of $A$ is equal to
 $2$ $\iff$ $\kappa=s\neq 0$ $\iff$ $(M,J_M,g)$ K\"ahler with $s\neq 0$ $\iff$
  $\mathcal I=\mathcal I_\frac{1}{2}=\{\mathbb J_\infty, \mathbb
 J_{-\infty} \}$ by Proposition~3. \qed

 \quad\\
 {\bf Proof of C~\&~D.} In those cases the matrix of $A$ in
a basis adapted to the decompostion $\C\otimes\bigwedge^+=\C
J_M\oplus^\perp(\bigwedge^{1,0}\oplus\bigwedge^{0,1})$ is
$\left[\begin{array}{ccc}x&0&0\\
0&y&0\\
0&0&y
\end{array}\right]$ with $y\neq0$. Let $f$ such that $\mathbb J_f\in
\mathcal I_{\frac{1}{2}}$,  $(m,Q)$ be any point of $Z$ and
$(\theta_1,...,\theta_4)$ be a local frame such that
$\left\{\begin{array}{l} J_M=\theta_1\wedge\theta_2+\theta_3\wedge\theta_4\\
Q\in Vect(I,J)\end{array}\right.$. So there exist $(a,b),
(\alpha,\beta,\gamma) \in\Sp^2$ such that $Q=aI+bJ$ and
$P=f(Q)=\alpha I+\beta J+\gamma K$. In that case at the point
$(m,Q)$ we have :
$$\begin{array}{lll}
 G(\theta_1,\theta_2)&=&0\\
 &=&[(1-\alpha^2)xI-\alpha\beta yJ-\alpha\gamma yK+Q(\gamma
 yJ-\beta yK ),Q]\\
 &=&\Big[\Big((1-\alpha^2)x-b\beta y\Big)I+(a-\alpha)\beta yJ+(a-\alpha)
 \gamma yK,Q]\\
 \end{array}
 $$
 $$
 \begin{array}{lll}
 \Longleftrightarrow
 \left\{
 \begin{array}{l}
 (a-\alpha)\gamma ya=0\\
 b\Big((1-\alpha^2)x-b\beta y\Big)=a(a-\alpha)\beta y
 \end{array}
\right. \\
\Longleftrightarrow
 \left\{
 \begin{array}{l}
 \gamma =0\\
\beta bx=y(1-a\alpha)\\
\alpha^2+\beta^2=1
 \end{array}
\right. &\textrm{ou} & \left\{
 \begin{array}{l}
 \alpha =a\\
 \beta=\frac{x}{y}b\\
 \beta^2+\gamma^2=b^2
 \end{array}
\right.\end{array}
$$
The resolution of $G(\theta_1,\theta_3)=0$ or
$G(\theta_1,\theta_4)=0$ gives the same system. Two cases can
happen first $\vert\frac{x}{y}\vert>1$ then the second system
doesn't have any solution and the first one has two solutions. An
easy computation enable us to verify that they correspond to
$f_1=u_1Id$ or $f_2=u_2Id$.

On the other hand if $\vert\frac{x}{y}\vert<1$ then the second
system gives two solutions which correspond to $f=e^{\pm
i\theta}Id$,
 whereas the first
system doesn't have any solution:
\[
\begin{array}{lll}
&1-\alpha^2=\beta^2
=\frac{y^2}{b^2x^2}(1-a\alpha)^2>\frac{(1-a\alpha)^2}{b^2}\\
\Longrightarrow&b^2-b^2\alpha^2>1+a^2\alpha^2-2a\alpha\\
\Longrightarrow&0>(\alpha-a)^2.
\end{array}
\]
When $\vert\frac{x}{y}\vert=1$ both system give the same
solutions. \qed

\subsection*{G) Example}
Let $\T$ be a torus which is a quotient of $\C$ by the
lattice $\Z\oplus i\Z$. Define $(M,g,I)$ to be the quotient
of the complex flat torus $\T^2=\T\times\T$ by the group
$H=\Z/2\Z$ generated by an element $h$. If
$(z_1,z_2)=(x_1+ix_2, x_3+ix_4)$ are the canonical
coordinates on $\C\times\C$, we have:
$$
h(z_1,z_2)=\Big(z_1+\frac{1}{2},-z_2\Big).
$$
The manifold $(M,g,I)$ is a bi-elliptic surface which is
scalar-flat K\"ahler; denote by $Z\lra M$ its twistor space. In
this section we will study in details this example, especially the
integrability of $\mathbb J_f$. Thanks  to Theorem~1, one knows
that $\mathbb J_{Id}$, $\mathbb J_\infty$ and $\mathbb J_{\lambda
Id}$ are integrable.

Let $(\frac{\partial}{\partial
  x_1},\frac{\partial}{\partial
  x_2},\frac{\partial}{\partial
  x_3},\frac{\partial}{\partial x_4})$ be the canonical
basis of $\C^2$ identified with $\R^4$. This furnishes a  basis of
vector fields on $\T^2$ and, consequently, a global trivialisation
of its twistor space $Z_0\simeq \T^2\times\Sp^2$. Define another
basis (on $\T^2$) by:
$$
  \theta_1+i\theta_2=\frac{\partial}{\partial
      x_1}+i\frac{\partial}{\partial x_2}\quad\textrm{ and }\quad
    \theta_3+i\theta_4=e^{2i\pi x_1}(\frac{\partial}{\partial
      x_3}+i\frac{\partial}{\partial x_4}).
$$
 Then,  $(\theta_1,\theta_2,\theta_3,\theta_4)$ is a global basis
on $\T^2$ which goes down to a basis of $M$. This
defines a new trivialisation of $Z_0$, denoted by
$\tilde M\times\Sp^2$.  The manifold $Z$ is the quotient of
$\tilde M\times \Sp^2$ by the group $\tilde H\simeq \Z/2\Z$,
generated by $\tilde h$ acting as follows:
\[
\begin{array}{cccc}
  \tilde h:&\tilde M\times\Sp^2&\lra&\tilde M\times\Sp^2\\
  &\Big(m,Q\Big)&\lms&\Big(h(m),Q\Big).
\end{array}
\]
Viewing $\Sp^2$ as a subspace of $\R\times\C$ with
coordinates $(a,z)$, the identity map $\Psi$ of $Z_0$ has
the following form in these trivialisations:
\[
\begin{array}{lclc}
  \Psi : & Z_0\simeq\T^2\times\Sp^2&\lra&Z_0\simeq\tilde
  M\times\Sp^2 \\
  &\xi\simeq\Big(m,(a,z)\Big)&\lms&\xi\simeq\Big(m,(a,e^{-2i\pi
    x_1}z)\Big).
\end{array}
\]
The matrix, in both basis $(\frac{\partial}{\partial
  x_1},\frac{\partial}{\partial
  x_2},\frac{\partial}{\partial
  x_3},\frac{\partial}{\partial x_4})$ and
$(\theta_1,\theta_2,\theta_3,\theta_4)$, of the natural complex
structure $I$ on $\T^2$  is equal to
$\left[\begin{array}{cccc} 0&-1&0&0\\
    1&0&0&0\\
    0&0&0&-1\\
    0&0&1&0
  \end{array}
\right] $.  According
to our notation, this is the  infinity section.\\
Endow $Z_0$ with the complex structure of twistor space $\mathbb
J_{Id}$. As $(\T^2,I)$ is hyperk\"ahler, the projection $
pr_2:Z_0\simeq \T^2\times\Sp^2\lra\C P^1$ is a holomorphic
submersion \cite{Boy88}. For $n \in \N^*$ and $\lambda \in \C^*$,
consider the application $f_n :\C P^1\lra \C P^1$ equal to
$\lambda z^n$. Then there exist two applications $f_1, f_2$
depending only on $|\lambda|$ such that:
\[
\xymatrix @R=0,2pc @C=0,5cm{
  \Sp^2&\ar[rr]^{f_n}&&&\Sp^2\\
  (a,z) \ar[dddd]&\ar[rr]&&&\Big(f_1(a),\lambda
  f_2(a)z^n\Big)
  \ar[dddd]\\
  \\
  \\
  \\
  \C\cup\{\infty\}&\ar[rr]^{f_n}&&&\C\cup\{\infty\}\\
  U=\frac{z}{1-a}&\ar[rr]&&&\lambda U^n }
\]
Introduce now the pull back $Z_n=f_n^\star Z_0$:
\[
\xymatrix @R=0,2pc @C=0,5cm{ Z_n \ar[dddd]&\ar[rr]&&&Z_0
  \ar[dddd]^{pr_2}\\
  \\
  \\
  \\
  \C P^1&\ar[rr]^{f_n}&&&\C P^1\\
}
\]
Since the fibration $Z_0\lra \C P^1$ is topologically
trivial, this is also the case for $Z_n\lra \C
P^1$. Therefore one can identify the manifold $Z_n$ with
$\T^2\times\Sp^2$ equipped with a complex structure denoted
by $J_n$. If one considers the morphism $\tilde f_n
=Id\times f_n : \T^2\times\Sp^2\lra \T^2\times\Sp^2$, which
respects the fibration, one has $J_n=\mathbb J_{\tilde
  f_n}$.

We were wondering whether this complex structure goes down
to $Z$, i.e.: does it commute with the action of the
group $\tilde H$?  We need to study $\Psi\circ \tilde
f_n\circ\Psi^{-1}$:
\[
\xymatrix @R=0,2pc @C=0,5cm{ \Big(\T^2\times\Sp^2,J_n\Big)
  &\ar[rr]^{\tilde f_n}&&&\Big(\T^2\times\Sp^2,\mathbb
  J_{Id}\Big)
  \\
  \Big(m,(a,e^{2i\pi x_1}z)\Big)\ar[dddd]_\Psi&\ar[rr]&&&\Big(m,\big(f_1(a),
  \lambda f_2(a)
  (e^{2i\pi x_1}z)^n\big)\Big)\ar[dddd]^\Psi\\
  \\
  \\
  \\
  \Big(\tilde M\times\Sp^2,J_n\Big)&\ar[rr]^{\Psi\circ\tilde
    f_n\circ\Psi^{-1}}&&&\Big(
  \tilde M\times\Sp^2,\mathbb J_{Id}\Big)\\
  \Big(m,(a,z)\Big)&\ar[rr]&&&\Big(m,\big(f_1(a),\lambda e^{2i\pi(n-1)
    x_1}f_2(a)z^n\big)\Big) }
\]
Thus, in the trivialisation of $Z_0\simeq \tilde M\times
\Sp^2$ associated to
$(\theta_1,\theta_2,\theta_3,\theta_4)$, the complex
structure $J_n$ is $\mathbb J_{\Psi\circ\tilde
  f_n\circ\Psi^{-1}}=\mathbb J_{\lambda e^{2i\pi(n-1) x_1}
  z^n}$. It commutes with $\tilde H$ if and only if $n$
is odd. Moreover, for $n$=1, $\tilde f_1$ is a
biholomorphism. We have proved the following:

\quad\\
{\bf Proposition 9.} {\it For all $\lambda\in\C^\star$ the complex
structures $\mathbb
  J_{\lambda z}$ on $Z$ are biholomorphic. Furthermore, the compatible almost
  complex structures $\mathbb J_{\lambda e^{2i\pi(n-1) x_1} z^n}$ are
  integrable for odd $n$. }

\quad\\
This proposition can be generalised to other bi-elliptic surfaces.
A computation similar to the one made in Proposition~5 enables us
to say that, for different integers $n$, these complex structures
are not deformation of each other. This is consequence of the fact
that they do not have the same Chern classes. Indeed, the first
Chern class satisfies $c_1(\mathbb J_{\lambda e^{2i\pi(n-1) x_1}
  z^n})=2(n+1)h$. In \cite{Des08}, following an idea of
LeBrun, we showed that for any hypercomplex manifold $M$ there
exist infinitely many complex structures on its twistor space
$Z\simeq M\times\Sp^2$ which are not deformation of each other.
Recall that the only compact hypercomplex surfaces are the torus,
the $K3$-surfaces and the quaternionic Hopf surfaces \cite{Boy88}.
The previous proposition can therefore be viewed as a
generalisation of this result to bi-elliptic surfaces.

\subsection*{H) Higher dimension}
The previous sections have focused on the $4$-dimensional case. We
now briefly give a generalization of Theorem~1 in higher
dimension.  Let $n>1$ and $(M,g)$ be an oriented $4n$-dimensional
Riemannian manifold, not necessarily compact. An almost
hypercomplex structure on $(M,g)$ is a triple $(I,J,K)$ of almost
complex structures compatible with the orientation and the metric,
such that $IJ=-JI=K$. When $I,J,K$ are integrable one speaks about
a hypercomplex structure. When they are K\"ahler one says that
$(M,g)$ is hyperk\"ahler.

An almost quaternionic structure $D$ on $(M,g)$ is a rank $3$
subbundle $D\subset \End(TM)$ which is locally spanned by an
almost hypercomplex structure $H=(I,J,K)$; such a triple is called
a local admissible basis. For $n>1$, one says that $(M,g,D)$ is a
quaternionic structure if there exists a torsion free connection
$\nabla$ on $TM$ preserving $D$. If one can choose $\nabla$ to be
the Levi-Civita connection, $(M,g,D)$ is called quaternionic
K\"ahler. This is equivalent to saying that the holonomy group of
$g$ is contained in $Sp(1)Sp(n)$ \cite{Bes87}.

A compatible almost complex structure on $(M,g,D)$ is a
section $J_M$ of $D\lra M$ such that $J_M^2=-Id$.

Let $(M,g,D)$ be a Riemannian almost quaternionic $4n$-manifold.
One can define a scalar product on $D$ by saying that a local
admissible basis of $D$ is orthonormal. One can then define the
twistor space $Z\lra M$, which is the unit sphere bundle of $D$.
This is a locally trivial bundle over $M$ with fiber $\Sp^2$ and
structure group $SO(3)$. As in the introduction, one can define a
natural metric $\tilde g$ and look for the compatible almost
complex structures on $(Z,\tilde g)$ which are integrable. When
$(M,g,D, J_M)$ is quaternionic K\"ahler with a compatible almost
complex structure $J_M$, its twistor space $(Z,\tilde g)$ admits
different compatible almost complex structures: $\mathbb
J_{\sigma}$, $\mathbb J_{Id}$, $\mathbb J_{\infty}$, $\mathbb
J_{\lambda Id}$, defined as previously. The main result of this
section is the following, where no hypothesis of compacity is
made.

\quad\\
{\bf Theorem 3.}  {\it Let $(M,g,D)$ be a quaternionic
  K\"ahler manifold.
\begin{itemize}
\item[{\rm A)}] The almost complex structure $\mathbb J_{\sigma}$
is
  never integrable.

\item[{\rm B)}]  The almost complex structure $\mathbb J_{Id}$ is
always
  integrable {\em \cite{Sal82}}.

 \item[{\rm C)}] If $(M,g,D,J_M)$ is a compatible almost complex
  quaternionic K\"ahler manifold the almost complex
  structure $\mathbb J_{\infty}$ is integrable if, and only
  if:
\begin{itemize}
\item[{\rm i)}]  $J_M$ is  integrable;

 \item[{\rm ii)}] $g$ is scalar-flat.
\end{itemize}

  \item[{\rm D)}] If $(M,g,D,J_M)$ is a quaternionic K\"ahler
  manifold with a compatible K\"ahlerian complex structure
  $J_M$ then, for all $\lambda\notin\{0,1\}$, the complex
  structure $\mathbb J_{\lambda Id}$ is integrable if, and
  only if, $g$ is scalar-flat.

  \item[{\rm E)}] Let $(M,g,D)$ be a quaternionic K\"ahler
  manifold. Then the scalar curvature is flat if, and only if,
  one (and then any) $m\in M$ has an open neighborhood
  $\mathcal U$ such that $(Z,\tilde g)$ admits over
  $\mathcal U$ an integrable compatible complex structure
  different from $\mathbb J_{Id}$.
  \end{itemize}}

\quad\\
Any quaternionic K\"ahler manifold which is scalar-flat is locally
hyperk\"ahler \cite{Bes87}. Thus, part~E of the previous theorem
yields a  characterization of locally hyperk\"ahler manifolds
among quaternionic K\"ahler's in terms of twistor spaces.

It is possible to give a simpler version of that theorem
in the compact case because of the following result.

\quad\\
{\bf Proposition \cite{Pon94}.} {\it In the compact case
  any compatible complex structure $J_M$ on a quaternionic
  K\"ahler manifold $(M,g,D)$ is automatically scalar-flat
  K\"ahler.}

\quad\\
In particular, in the compact case, Theorem~3 has the following
corollary.

\quad\\
{\bf Corollary 2.} {\it Let $(M,g,D,J_M)$ be a compact
  quaternionic K\"ahler manifold with a compatible almost
  complex structure. Then $J_M$ is integrable if, and only if,
  $\mathbb J_\infty$ is integrable. In this case $\mathbb
  J_{\lambda Id}$ is integrable for all $\lambda\in\C^\star$.
}

\quad\\
{\bf Proof of Theorem~3.} Proposition~1 and Proposition~2 remain
true in dimension $4n$.  Since $\sigma$ is an antiholomorphic
involution when restricted to the fibers, part A can be  easily
proved.

The proof of part~B is the same as the one given in dimension $4$.
Notice first that $d\pi F_{ij}=-E(\theta_i,\theta_j)$ for all
$(i,j)\in\{1,\dots,4n\}$. It remains to show that
$G(\theta_i,\theta_j)=0$ for all $i,j\in\{1,\dots,4n\}$. To get
that result we use the following lemma.

\quad\\
{\bf Lemma 3 \cite{Bes87}.}  {\it Let $r(.,.)$ be the Ricci
  tensor. For all $(X,Y)\in TM$ one has:
  \[
  \begin{array}{l}
    \,[R(X,Y),I]=\gamma(X,Y)J-\beta(X,Y) K\\
    \,[R(X,Y),J]=-\gamma(X,Y)I+\alpha(X,Y)K\\
    \,[R(X,Y),K]=\beta(X,Y)I-\alpha(X,Y)J\\
  \end{array}
  \textrm{ with } \left\{
    \begin{array}{l}
      \alpha(X,Y)=\frac{2}{n+2}r(IX,X)\\
      \beta(X,Y)=\frac{2}{n+2}r(JX,X)\\
      \gamma(X,Y)=\frac{2}{n+2}r(KX,X)
    \end{array}
  \right.
  \]}
\quad\\
 Let $(m,I)\in Z$ and $(I,J,K)$ be a local admissible
basis. Then Lemma~3 yields:
\[
\begin{array}{lll}
  G(\theta_i,\theta_j)&=&\Big[R\Big(\theta_i\wedge
  \theta_j-I\theta_i\wedge I\theta_j\Big)+IR\Big(\theta_i\wedge
  I\theta_j+I\theta_i\wedge\theta_j\Big),I\Big]\\
  &=&\gamma(\theta_i,\theta_j)J-\beta(\theta_i,\theta_j) K-
  \gamma(I\theta_i,I\theta_j)J+\beta(I\theta_i,I\theta_j) K \\
  &&+ \gamma(I\theta_i,\theta_j)K+\beta(I\theta_i,\theta_j) J+
  \gamma(\theta_i,I\theta_j)K+\beta(\theta_i,I\theta_j) J\\
\end{array}
\]
But any quaternionic K\"ahler manifold is Einstein \cite{Ber66},
hence $r=\frac{s}{4}g$, where $s$ is the scalar curvature of $g$.
One then has, for all $(\theta_i,\theta_j)$:
\[
\begin{array}{lll}
  G(\theta_i,\theta_j)&=&\frac{2s}{4(n+2)}
\Big(\big(2g(K\theta_i,\theta_j)-2g(K\theta_i,\theta_j)\big)J+
  \big(2g(J\theta_i,\theta_j)-2g(J\theta_i,\theta_j)\big)K\Big)
  \\
  &=&0.
\end{array}
\]

To prove part~C observe that, as in dimension $4$: \{$\mathbb
J_\infty$ integrable\} $\iff$
\{$E(\theta_i,\theta_j)=G(\theta_i,\theta_j)=0$\} $\iff$ \{$J_M$
integrable and $G(\theta_i,\theta_j)=0$\}. Since $(M,g,Q)$ is
Einstein, $(M,g,Q)$ scalar-flat implies $(M,g,Q)$ Ricci-flat and
$G(\theta_i,\theta_j)=0$. The converse is a consequence of part~E:
if $\mathbb J_\infty$ integrable then $s=0$.

To get part~D we use the technique (and notation) of dimension
$4$. Let $z\in\pi^{-1}(m)$ be a point in $Z$ over $m\in M$. As
$J_M$ is a parallel compatible complex structure, one can build a
basis $(\theta_1,\dots,\theta_{4n})$ over an open set $\mathcal U$
of $M$ such that $\nabla\vert_m\theta_i=0$ and the matrix of $J_M$
in that basis is $I$. This frame gives a local admissible basis
$(I,J,K)$ of $D$ and therefore a local trivialisation of $Z$ over
$\mathcal U$, in which the restrictions of $f$ to the fibers do
not depend on the second variable. So $F_{ij}\vert
_m=E(\theta_i,\theta_j)\vert_m=0$ and $G(\theta_i,\theta_j)=0$ if
$s$=0. The converse is again a consequence of part~E.

Proof of~E: suppose that the scalar curvature $s$ of
$(M,g,D)$ is non zero. Let $f:Z\lra Z$ be a morphism such
that $\mathbb J_f$ is integrable over an open set $\mathcal
U$. Let $(m,Q)$ be a point in $\pi^{-1}(\mathcal U)$ and set
$f(m,Q)=P$.  If $\mathcal U$ is small enough there
exists an orthonormal basis $(\theta_1,\dots,\theta_{4n})$ and
a local admissible basis $(I,J,K)$ such that $P=J$. Write
$Q=aI+bJ+cK$ with $(a,b,c)\in\Sp^2$.

As $\mathbb J_f$ is integrable we have
$G(\theta_1,\theta_2)=0$ everywhere. In particular at the
point $(m,Q)$ :
\[
\begin{array}{lll}
  G(\theta_1,\theta_2)&=&0\\
  &=&\Big[R(\theta_1\wedge\theta_2+\theta_3\wedge\theta_4)-Q
  R(\theta_1\wedge\theta_4+\theta_2\wedge\theta_3),Q\Big]\\
  &=&\frac{2}{n+2}(-2cJ+2bK)-Q\Big[R(\theta_1\wedge\theta_4+\theta_2\wedge\theta_3),
  Q\Big]  \\
  &=&\frac{4}{n+2}\Big(-cJ+bK-Q(-bI+aJ)\Big)\\
  &=&\frac{4}{n+2}\Big(acI+c(b-1)J+(b-1)K\Big)
\end{array}
\]
Hence $Q=J=P$ for any $(m,Q)\in \pi^{-1}(\mathcal U)$, that is
$f=Id$.

The converse is the same as the one given in section~E. Indeed, a
quaternionic K\"ahler manifolds $(M,g,D)$ admits, locally,
infinitly many compatible complex structures $J_M$ (for example
\cite{AMP98}). \qed

\bibliographystyle{alpha}
\bibliography{bibstar}

\end{document}